\documentclass[10pt,]{article}
\usepackage{xspace}
\usepackage{amssymb}
\usepackage{amsmath,amscd}
\usepackage{amsfonts}
\usepackage{amsthm}
\usepackage{keyval} 
\usepackage{array}
\usepackage[all,cmtip]{xy}
\usepackage{setspace}
\usepackage{enumerate}

\newtheorem{theorem}{Theorem}[section]
\newtheorem{definition}{Definition}[section]

\begin{document}
\title{Categories of modules and their deformations}
\author{Romie Banerjee}
\date{}

\maketitle

\begin{abstract}
We develop an obstruction
theory for lifting compact objects to the stable $\infty$ category of
quasi-coherent modules over a derived geometric stack $X$ from the
category of modules over its underlying classical stack $X^{cl}$. The
obstructions live in Andre-Quillen cohomology. 
\end{abstract}

\section{Introduction}
The derived category of quasi-coherent modules over a scheme or an algebraic stack is usually very badly behaved in the sense that it is not controlled by a small data. In certain cases it is possible to find a set of compact generators of the derived category of modules in question. For example, if $R$ is a commutative ring then the triangulated category $D(R)$ of chain complexes of $R$-modules modulo weak equivalences of chain cohomology isomorphisms is compactly generated. Similar thing is true of the unbounded derived category of quasi coherent modules over a quasi-compact separated scheme \cite{Neeman}. In general not all algebraic stacks have this property. Stable homotopy theory gives rises to more sources of interesting trianglated categories. For any $E_{\infty}$-ring spectrum $A$ the derived category of $A$-modules is compactly generated. For any derived scheme, formed by gluing derived affine schemes $\hbox{Spec}(A)$ along Zariski maps of $E_{\infty}$-rings, the derived category of quasi-coherent modules form a compactly generated triangulated category \cite{Drinfeld}. We are interested in the triangulated category of the derived category of quasi-coherent modules over any general derived algebraic stack. Throughout this paper we think of derived algebraic stacks, once rigidified, as being equivalent to cosimplicial connective $E_{\infty}$-rings.

Given a derived $\infty$-stack $X$, we want to study the stable $\infty$-category of quasi-coherent modules over $X$. If $X$ is an algebraic stack, i.e. $X$ admits an atlas by simplicial derived affive scheme $U_{\bullet}$, we get a cosimplicial stable $\infty$-category $\hbox{Mod}(U_{\bullet})$. The stable $\infty$-category modules over the stack $X$ is the totalization $\hbox{Tot}(\hbox{Mod}(U_{\bullet}))$.
 
Let $\mathfrak{C}_{\infty}$ be the category of connective $E_{\infty}$ rings. Another way to approach this is to consider the stack $QC$ considered as a moduli functor $$QC: \mathfrak{C}_{\infty} \rightarrow Pr_{st-\infty}$$ where the right side is the $\infty$-category of presentable stable $\infty$-categories, so that $QC(A) = \hbox{Mod}(A)$ and $QC$ takes a map of modules $f:A\rightarrow B$ to the functor $-\otimes_AB$. This naturally extends to a functor between $\infty$-categories. The desired object, i.e. the $\infty$-category of quasi-coherent modules over any $\infty$-stack $X:\mathfrak{C}_{\infty} \rightarrow SSet$ is $\hbox{Hom}_{\infty\hbox{-stacks}}(X,QC)$, the hom space in the $\infty$-topos of $\infty$-stacks.

If $A$ is a connective $E_{\infty}$ ring which admits a postnikov tower decomposition and $\mathcal{M}$ an $\infty$-stack which admits a cotangent complex and is infinitesimaly cohesive \cite{HAG}, lifting a family of objects classified by $\mathcal{M}$ on the ordinary affive scheme $\hbox{Spec}\,\pi_0A$ to the derived affine scheme $\hbox{Spec} A$ is a problem in deformation theory. It is controlled by the cotangent complex of the stack $\mathcal{M}$. Associated to any derived algebraic $\infty$-stack $X$ there is an ordinary algebraic $\infty$-stack $X^{cl}$ which admits an atlas of a cosimplicial ordinary ring obtained by taking sectionwise $\pi_0$ of the the atlas of $X$. We can think of $X$ as an infinitesimal extension of the underlying $X^{cl}$. 

Let $X$ be a derived $\infty$-stack. Let $X^{cl}$ be it's associated classical (non-derived stack). There is a natural map $i:X^{cl} \rightarrow X$. The induced map on derived categories: $$D_{qc}(X) \rightarrow D_{qc}(X^{cl})$$

Given arbitrary $x$ in $D_{qc}(X)$ and a perfect $u$ in $D_{qc}(X^{cl})$, with a map $u \rightarrow x$ we would like to find cohomological obstructions for lifting $u$ to a perfect module $\widetilde{u}$ over $X$ and a map $\tilde{u} \rightarrow x$ over $X$ which restricts to $u \rightarrow x$ over $X^{cl}$

The main result is

\begin{theorem}
Let $X$ be a perfect derived algebraic $n$-stack for some $n$ and let $\widetilde{X}$ be a square-zero extension of $X$. Let $x:\widetilde{X} \rightarrow QC$ be a complex of quasi-coherent modules over $\widetilde{X}$ and let $u:X \rightarrow QC^{\omega}$ be a complex of perfect modules over $X$, along with a map $u \rightarrow x$ in $QC(X)$.

\begin{itemize} 

\item Then there exists an obstruction theory for deforming $u$ to a $\widetilde{u}:X \rightarrow QC^{\omega}$. The space of deformations is isomorphic to $\Omega\hbox{Hom}_{\mathcal{O}_X}(\alpha^*L_{QC^{\omega}},N)$ with loops based at the trivial derivation. 

\item If this space in non-empty and $\widetilde{u}$ is a deformation of $u$, then there exists a perfect module $y_{\beta}: X \rightarrow QC^{\omega}$ along with maps $\beta:u \rightarrow y_{\beta}$ and $y_{\beta} \rightarrow x$ in $QC(X)$ such that the triangle commutes in $QC(X)$

$$\xymatrix{
u \ar[dr]_{\beta} \ar[rr] &&x \\
&y_{\beta} \ar[ur]
}$$

There is an obstrucion theory for lifting $\beta$ to $\widetilde{\beta}:\widetilde{u} \rightarrow \widetilde{y}_{\beta}$ such that $\widetilde{u} \rightarrow \widetilde{y}_{\beta} \rightarrow x$ is a deformation of $\alpha:u \rightarrow x$. 

More precisely, there exists a moduli functor $\mathcal{G}:\Omega_{u,y_{\beta}}QC_{/X\times X}$ and an cocycle in the Andre-Quillen cohomology 

$$\alpha(u,y_{\beta}) \in \hbox{Hom}_{\mathcal{O}_X}(\beta^*L_{\mathcal{G}},N)$$ such that, if $\alpha(u,y_{\beta})=0$ there exists a lift $\widetilde{\beta}$. The space of all such deformations is isomorphic to $$\Omega\hbox{Hom}_{\mathcal{O}_X}(\beta^*L_{\mathcal{G}},N)$$ where the loops are based at the trivial derivation.
\end{itemize}

\end{theorem}

\tableofcontents
\section{Derived $\infty$-stacks: Overview}
In this section we give a brief introduction to geometric $\infty$-stacks. The reader may find the necessary details on $ \infty$-categories an $\infty$-topoi in \cite{HTT}.

Let $\mathcal{C}^{op}$ denote a presentable $\infty$-category (connective $E_{\infty}$ rings, or simplicial commutative rings) with an topology $\tau$ on $\mathcal{C}$. A derived $\infty$-prestack is a functor
$$\mathcal{F}: \mathcal{C}^{op} \rightarrow \mathcal{S}$$
$$\hbox{Spec}\,A = \hbox{Hom}_{\mathcal{C}^{op}}(A,-)$$

$\mathcal{F}$ is an $\infty$-stack if it satifies Cech descent with respect to $\tau$; $X \in \hbox{Fun}^L(\mathcal{C}^{op},\mathcal{S})$ and $\mathcal{F}$ takes the Cech nerve of any $\tau$-cover $U\rightarrow X$ to a limit diagram.

$\mathcal{F}$ is an algebraic $\infty$-stack if there is a cosimplicial object $A_{\bullet} \in (\mathcal{C}^{op})^{N(\Delta)}$ and $\mathcal{F}(B) = |\hbox{Hom}_{\mathcal{C}^{op}}(A_{\bullet},B)|$, $$\mathcal{F} = \hbox{colim}_{\Delta^{op}}\, \hbox{Spec}\,A_{\bullet}$$ in the $\infty$-category of $\infty$-stacks.

\subsection{The quasi-coherent $\infty$-stack} 

A quasi-coherent sheaf on a scheme $X$ is a morphism of stacks $X \rightarrow Mod$ from $X$, considered as a stack, into the canonical stack $$Mod: \hbox{Spec}A \mapsto \hbox{Mod}_A$$ of modules which corresponds to the bifibration $$T_{CRings} \simeq Mod \rightarrow CRings$$ from the tangent category of the the category of commutative rings to commutative rings.

This definition of quasi-coherent sheaves generalizes to any $(\infty,1)$-topos, and over arbitrary $\infty$- sites. Let $\mathcal{C}$ be symmetric monoidal $\infty$-category equipped with Grothendieck $\infty$ topology such that $\mathcal{C}^{op}$ is presentable. The tangent $\infty$ category $T(\mathcal{C}^{op}) \rightarrow \mathcal{C}^{op}$ is the bifibration whose fibers over an object $A \in \mathcal{C}$ plays the role of the $\infty$-groupoid of modules over $A$, see \cite{DAGIV}.

Under the $\infty$-Grothendieck construction this corresponds to a $(\infty,1)$ presheaf $$Mod_{\infty}: \mathcal{C}^{op} \rightarrow \widehat{\hbox{Cat}_{\infty}}$$ $$Mod_{\infty}: \hbox{Spec} R \mapsto \hbox{Stab}(\mathcal{C}^{op}_{/R})$$ where $\hbox{Spec} R$ for $R \in \mathcal{C}^{op}$ is the affine object in the geometry defined over $\mathcal{C}^{op}$, or directly in terms of test spaces $$Mod_{\infty}: \mathcal{U} \mapsto \hbox{Stab}(\mathcal{C}^{\mathcal{U}/})$$ This makes $Mod_{\infty}^{\mathcal{C}} \in \hbox{Shv}_{(\infty,1)}^{\widehat{\hbox{Cat}_{\infty}}}(\mathcal{C}) = [\mathcal{C}^{op}, \widehat{\hbox{Cat}_{\infty}}]$.

Let $H = \hbox{Shv}_{\infty}(\mathcal{C})$ be the $\infty$-topos of $\infty$-stacks on $\mathcal{C}$ and $X \in H$ be an $\infty$-stack. The stable $\infty$ category of quasi coherent modules over $X$ is the Hom space in the $\infty$-topos $H$; 
\begin{definition}
\begin{equation}
QC(X) = \hbox{Hom}_H(X, Mod_{\infty})
\end{equation}
\end{definition}

Notice that $H \subset [\mathcal{C}^{op}, \widehat{\hbox{Cat}_{\infty}}]$ as any $\infty$-groupoid is in $\widehat{\hbox{Cat}_{\infty}}$.

$QC(X)$ is computed using the Yoneda-Kan extension. 

$$\xymatrix{ 
\mathcal{C} \ar[r]^j \ar[d]_{Mod} &P(\mathcal{C}) \ar[dl]^{\hbox{Kan}_j(Mod)} \\
\widehat{\hbox{Cat}}
}$$ 

$$\xymatrix{ 
\hbox{Fun}(\mathcal{C}^{op},\mathcal{D}) \ar@<1ex>[r]^{\hbox{Kan}}  &\hbox{Fun}^R (P (\mathcal{C})^{op}, \mathcal{D}) \ar@<1ex>[l]^j 
}$$

By definition $\hbox{Kan}(F)(Y) = \hbox{lim}_{j(U) \rightarrow Y} F(U)$, where $Y = \hbox{colim}_{j(U)\rightarrow Y} j(U)$ in $P(\mathcal{C})$. The above adjunction is an equivalence of $\infty$ categories; it follows from the the standard adjunction $\xymatrix{ \hbox{Fun}(\mathcal{A}, \mathcal{B}) \ar@<1ex>[r]^{\hbox{Kan}} &\hbox{Fun}^L(P(\mathcal{A}),\mathcal{B}) \ar@<1ex>[l]^j}$ being an equivalence of $\infty$ categories.

For a prestack $X \in P(\mathcal{C}) = \hbox{Fun}(\mathcal{C}^{op}, \hbox{SSets})$, suppose $X = \hbox{colim}_{\alpha} j(\hbox{Spec‭}R_{\alpha})$ the $$\widetilde{QC}(X) = \hbox{Kan}_j(Mod)(X) = \hbox{lim}_{\alpha} Mod (\hbox{Spec}R_{\alpha}) = \hbox{lim}_{\alpha} \hbox{Stab}(\mathcal{C}^{op}_{/R_{\alpha}}).$$

If $X$ is an $\infty$-stack, $QC(X)$ can be expressed as a limit similarly, $$QC: \xymatrix{\infty{\hbox{-stacks}} \ar[r]^{i^{op}} &P(C)^{op} \ar[r]^{\widetilde{QC}} &&\widehat{\hbox{Cat}} }$$ however since $i^{op}$ doesn't preserve limits, it is not straightforward to show.

If $X$ is a geometric $\infty$-stack (i.e. atlas by a simplicial object in $\mathcal{C}$), we want to compute $QC(X)$. $QC$ is the composition $$QC: \xymatrix{\hbox{geometric-}\infty\hbox{-stacks}^{op} \cong \mathcal{C}^{\Delta} \ar[rr]^{i^{op}} &&P(\mathcal{C}^{op}) \ar[r]^{\widetilde{QC}} \ar[r]^{\widetilde{QC}} &&& \widehat{\hbox{Cat}_{\infty}} }$$

If $A_{\bullet} \in \mathcal{C}^{\Delta}$ is the cosimplicial object such that simplicial object in $\mathcal{C}$, $\hbox{Spec}(A_{\bullet})$ (or simply, the simplicial affine scheme) is an atlas for $X$, then 
$$i(A_{\bullet}) = \hbox{Hom}_{\mathcal{C}^{op}}(A_{\bullet},-)$$ 

That is, as an object in the prestack category $i(A_{\bullet})$ evaulates on objects in $\mathcal{C}^{op}$ as the geometric realization $$ i(A_{\bullet})(R) = |\hbox{Hom}_{\mathcal{C}^{op}}(A_{\bullet},R)|.$$

or, $i(A_{\bullet}) = \hbox{colim}_{\Delta^{op}}\hbox{Spec}(A_{\bullet})$ in the $\infty$ category of affine $\mathcal{C} $-schemes.

Therefore, 
\begin{equation}
QC(A_{\bullet}) = \widetilde{QC}(i^{op}(A_{\bullet})) = \widetilde{QC}(\lim_{\Delta}\hbox{Spec}(A_{\bullet})) = \lim_{\Delta}\widetilde{QC}(\hbox{Spec}A_{\bullet}) = \hbox{Tot} \hbox{Mod}(A_{\bullet})
\end{equation}
where the limit/Tot is taken in the category of the stable presentable $\infty$ categories. 

\section{Deformation Theory} 

In this section we describe the basic setup for doing deformation theory of geometric $\infty$-stacks. We will closely follow Lurie's DAG IV \cite{DAGIV}.

Let $\mathcal{D}$ be a presentable $\infty$ category, then the tangent category $T_{\mathcal{D}}$ is the fiberwise stablization of the projection map $$\hbox{Fun}(\Delta^1, \mathcal{D}) \rightarrow \hbox{Fun}(\{1\},\mathcal{D}) \simeq \mathcal{D}$$ 
Roughly speaking, an object of the tangent bundle $T_{\mathcal{D}}$ consists of a pair $(A,M)$, where $A \in \mathcal{D}$ and $M \in \hbox{Stab}(\mathcal{D}_{/A})$; here Stab is the stabilization construction applied to an $\infty$ category. If $\mathcal{D}$ is the ordinary category of commutative rings(replace stabilization with abelianization) then the associated tangent category is equivalent to the category of modules; the objects are pairs $(A,M)$, where $A$ is a commutative ring and $M$ is a $A$-module. If $\mathcal{D}$ is the $\infty$-category of $E_{\infty}$-rings or simplicial commutative rings then the tangent category recovers the categories of modules over such objects. Using this analogy, we can define a {\it module} over an object $A$ to be an object of the fiber of the tangent category $T_{\mathcal{D}}$ over $\mathcal{D}$, ie. the stable $\infty$-category $T_{\mathcal{D}} \times _{\mathcal{D}} A \simeq \hbox{Stab}(\mathcal{D}_{/A})$.

The {\it cotangent complex functor} $L:\mathcal{D} \rightarrow T_{\mathcal{D}}$ is the left adjoint to the forgetful functor
$$T_{\mathcal{D}} \rightarrow \hbox{Fun}(\Delta^1, \mathcal{D}) \rightarrow \hbox{Fun}(\{0 \}, \mathcal{D}) \simeq \mathcal{D}$$ such that the cotangent complex $L_A$ of and object $A$ is in $\hbox{Stab}(\mathcal{C}_{/A})$. In other words, the composition $$ \mathcal{C} \rightarrow ^L T_{\mathcal{D}} \rightarrow \hbox{Fun}(\Delta^1, \mathcal{D}) \rightarrow \hbox{Fun}({1}, \mathcal{D}) \simeq \mathcal{D}$$ is the identity functor. 

The absolute cotangent complex functor $L: \mathcal{D} \rightarrow T_{\mathcal{D}}$ is defined to be the compostion 
$$\mathcal{D} \rightarrow \hbox{Fun}(\Delta^1, \mathcal{D}) \rightarrow T_{\mathcal{D}}$$ 
where the first map is the given by the the diagonal embedding and the second map is the left adjoint to the forgetful functor $G:T_{\mathcal{D}} \rightarrow \hbox{Fun}(\Delta^1, \mathcal{D})$

$$\xymatrix{
T_{\mathcal{D}} \ar[rr]^G \ar[dr]_p &&\hbox{Fun}(\Delta^1, \mathcal{D}) \ar[dl]^{\hbox{ev}_1}\\
&\mathcal{D}
}$$

Since the diagonal embedding is left adjoint to the evaluation map $\hbox{Fun}(\Delta^1, \mathcal{D}) \rightarrow \hbox{hbox}(\{ 0 \}, \mathcal{D})$ the absolute cotangent complex functor is left adjoit to the composition $T_{\mathcal{D}} \rightarrow \hbox{Fun}(\Delta^1, \mathcal{D}) \rightarrow \hbox{Fun}(\{0\}, \mathcal{D})$.

The fiber of the tangent bundle $T_{\mathcal{D}}$ over $A \in \mathcal{D}$ can be identified with the stable envelope $\hbox{Stab}(\mathcal{D}_{/A})$. Under this identification the cotangent complex $L_A \in \hbox{Stab}(\mathcal{D}_{/A})$ corresponds to the image of $\hbox{id}_A \in \mathcal{D}_{/A}$ under the suspension functor $$ \Sigma^{\infty}: \mathcal{D}_{/A} \rightarrow \hbox{Stab}(\mathcal{D}_{/A}).$$
 
The {\it trivial square zero extension} of $A \in \mathcal{D}$ along a $A$-module $M$, denoted by $A\oplus M$ is the image of the $M$ under the functor $$\Omega^{\infty}:\hbox{Stab}(\mathcal{D}_{/A}) \rightarrow \mathcal{D}_{/A} \rightarrow \mathcal{D}$$

Given an object $A \in \mathcal{D}$ and a $A$-module $M \in T_{\mathcal{D}} \times_{\mathcal{D}} \{A \}$, a {\it derivation} of $A$ into $M$ is a map $\eta: L_A \rightarrow M$ in the $\infty$-category $T_{\mathcal{D}} \times_{\mathcal{D}} \{ A \}$. The derivation $\eta$ equivalently gives a map from $A$ to the trivial square-zero extension of $A$ defined by $M$ in the category $\mathcal{D}$, $$d_{\eta}: A \rightarrow A\oplus M$$ The derivation classified by the zero map $L_A \rightarrow M$ (this is a stable category) corresponds a canonical section $d_0:A \rightarrow A\oplus M$ in $\mathcal{D}$. The {\it square-zero extension} of $A$ defined by the derivation $\eta: L_A \rightarrow M$ is the pullback in the $\infty$ category $\mathcal{D}$. 

$$\xymatrix{
A^{\eta} \ar[r] \ar[d] &A \ar[d]^{d_0} \\
A \ar[r]^{d_{\eta}} &A\oplus M
}$$

Let $f:\tilde{A} \rightarrow A$ be a morphism in $\mathcal{D}$. Then $f$ is a square-zero extension if there exists a derivation $\eta:L_A \rightarrow M$ and an equivalence $\tilde{A} \simeq A^{\eta}$ in the $\infty$-category $\mathcal{D}_{/A}$. 

The square-zero extension $\tilde{A}$ will also be alternatively denoted by $A\oplus_{\eta}\Omega M$, so that $A \oplus_0 \Omega M \simeq A \oplus M$.

\subsection{Infitesimal Extensions of $\infty$-stacks}

Suppose $A_{\bullet}$ is a cosimplicial object in $\mathcal{C}^{op}$ and $X = \hbox{colim}_{\Delta^{op}} \hbox{Spec}A_{\bullet}$ the associated algebraic stack. Then $$\hbox{Mod}_{\mathcal{O}_X} = \hbox{Hom}_{H}(X, Mod)$$ We have seen have how to compute this 

$$\hbox{Mod}_{A_{\bullet}} = \hbox{Tot}_{[n] \in \Delta}\hbox{Stab}(\mathcal{C}^{op}_{/A_n}) \simeq \hbox{Tot}_{[n] \in \Delta} \hbox{Mod}_{A_n}$$

Therefore a module over $\mathcal{O}_X$ is an object in the totalization of a cosimplicial stable $\infty$-category. The $0$-simplices of the Tot stable $\infty$-category are exactly $A_0$-modules $+$ descent data, i.e. $\mathcal{O}_X$-modules. A $\mathcal{O}_X$-module $N$ is a cosimplicial diagram of modules, $N_{n} \in \hbox{Stab}(\mathcal{C}^{op}_{/A_n})$, and descent data. The trivial square zero extension defined by each $A_n$-module $N_n$ is the image of $N_n$ under the map $\hbox{ev}_0\circ \Omega^{\infty}: \hbox{Stab}(\mathcal{C}^{op}_{/A_n}) \rightarrow \mathcal{C}^{op}$. 

Let $\hbox{Stab}(\mathcal{C}^{op}_{/A_{\bullet}})$ denote the cosimplicial stable $\infty$-category induced by the cosimplicial diagram $A_{\bullet}$; given a map $A \rightarrow B$ in $\mathcal{C}^{op}$ there is a naturally induced map of stable $\infty$-categories $$ \hbox{Stab}(\mathcal{C}^{op}_{/A}) \rightarrow \hbox{Stab}(\mathcal{C}^{op}_{/B}).$$ 

The limit of this diagram in the $\infty$-category of stable $\infty$-categories is the category whose objects consists of an object in each category and descent data required to glue them. In the general situation we will use the notation $\mathcal{D}_{/A\rightarrow B}$ for the $\infty$-category $\hbox{lim}(\mathcal{C}_{/A} \rightarrow \mathcal{C}_{/B})$. Therefore in our case of interest, the $\infty$-category $\hbox{Stab}(\mathcal{C}^{op}_{/A_{\bullet}})$ (where Stab is taken level-wise) is the limit category $\hbox{Tot}_{[n] \in \Delta}\hbox{Stab}(\mathcal{C}^{op}_{/A_[n]})$.

We can apply the functor $\Omega^{\infty}$ to the cosimplicial stable presentable $\infty$-category and compose with evaluation at $\{ 0 \} \in \Delta^1$. 
$$\xymatrix{
\hbox{Stab}(\mathcal{C}^{op}_{/A_{\bullet}}) \ar[r]^{\Omega ^{\infty}} &\mathcal{C}^{op}_{A_{\bullet}} \ar[r]^{\hbox{ev}_0} &\mathcal{C}^{\Delta^{op}}
}$$

Let $A_{\bullet} \in (\mathcal{C}^{op})^{\Delta}$ and $N$ a $A$-module, that is an object in the totalization of the cosimplicial category $\hbox{Stab}(\mathcal{C}^{op}_{/A})$. The the {\it trivial square-zero extension} of $A_{\bullet}$ defined by $N$ is the image of $N$ under the map $\Omega^{\infty}\circ \hbox{ev}_0$. Denote this cosimplicial object in $\mathcal{C}^{op}$ by $A_{\bullet} \oplus N$. 

If $X$ is the geometric $\infty$-stack whose atlas is the simplcial affine $\mathcal{C}$-scheme $\hbox{Spec}A_{\bullet}$, we'll denote the trivial square zero extension by $\mathcal{O}_X\oplus N$.

The absolute cotangent complex of a cosimplicial ring $A_{\bullet}$ is the absolute cotangent complex of the associated geometric stack $X = \hbox{colim}_{\Delta^{op}} \hbox{Spec}A_{\bullet}$, $L_X$ (defined in the next section).

$L_X \in \hbox{Stab}(\mathcal{C}^{op}_{/A_{\bullet}})$. For any $\mathcal{O}_X$-module $N$, a {\it derivation} of $X$ into $N$ is a map on the stable$\infty$-category $\hbox{Stab}(\mathcal{C}^{op}_{/A_\bullet})$ $$\eta: L_X \rightarrow N$$ By adjunction, this is equivalent to giving a map $A_{\bullet} \rightarrow A_{\bullet}\oplus N$ in $\mathcal{C}^{\Delta^{op}}$. 

The {\it square-zero extension of $A_{\bullet}$ defined by $\eta$} is the pullback in the $\infty$-category $\mathcal{C}^{\Delta^{op}}$ 
$$\xymatrix{
A_{\bullet}^{\eta} \ar[r] \ar[d] &A_{\bullet} \ar[d]^{d_0} \\
A_{\bullet} \ar[r]^{d_{\eta}} & A_{\bullet}\oplus N
}$$

Denote the geometric $\infty$-stack defined by the atlas $\hbox{Spec}A_{\bullet}^{\eta}$ by $X\oplus_{\eta}[\Omega N]$.

\subsection{Cotangent complexes of $\infty$-stacks}

The {\it cotangent complex of an $\infty$-stack}. Let $F$ be an $\infty \mathcal{C}$-stack, i.e. an object in $\hbox{Fun}(\mathcal{C}^{op}, \hbox{SSet})$. For $A \in \mathcal{C}^{op}$ and $M \in \hbox{Stab}(\mathcal{C}^{op}_{/A})$. Let $A\oplus M$ be the trivial square-zero extension of $A$ by $M$. Let $$x:\hbox{Spec}A \rightarrow F$$ be a $A$-point. Fix the following notation 
$$X:= \hbox{Spec}A$$ $$X[M]:= \hbox{Spec}(A\oplus M)$$ The natural augmentation $A\rightarrow A\oplus M$ gives a natural map of stacks $X \rightarrow X[M]$. 

The {\it space of derivarions} from $F$ to $M$ at $x$ is defined by $$\hbox{Def}_F(x,M) := \hbox{Hom}_{X/Aff_{\mathcal{C}}}(X[M],F)$$

As $M \mapsto X[M]$ is functorial in $M$ is functorial in $M$, there is a well defined functor $$\hbox{Def}_F(x,-): \hbox{Mod}_A \rightarrow \hbox{SSets}$$

defined to be the homotopy fiber in the $\infty$-category of simplicial sets
$$\xymatrix{
\hbox{Def}_F(x,M) \ar[r] \ar[d] &F(X[M]) \ar[d]\\
\star \ar[r]^x &F(X)
}$$

The map $X \rightarrow X[M]$ has a canonical section ( the zero derivation $d_0:A \rightarrow A \oplus M$). Therefore, $\hbox{Def}_F(x,M)$ is a pointed space.

$F$ has a cotangent complex at $x$ if the functor $\hbox{Def}_F(x,M)$ is corepresented by a $A$-module $L_{F,x}$. The module $L_{F,x} \in \hbox{Stab}(\mathcal{C}^{op}_{/A})$ is the cotangent complex of $F$ at $x$. 

The $\infty$ stack $F$ has an {\it absolute cotangent complex} if for any $A \in \mathcal{C}^{op}$ and any $x \in F(A)$, $F$ has a cotangent complex $L_{F,x}$ at $x$ and for any commutative diagram in $\infty\hbox{-stacks}_{/F}$ 

$$\xymatrix{
\hbox{Spec}A \ar[rr]^u \ar[dr]_x && \hbox{Spec}B \ar[dl]^{x'} \\
&F
}$$
the natural morphism $u^*L_{F,x'} \rightarrow L_{F,x}$ is an equivalence in $\hbox{Stab}(\mathcal{C}^{op}_{/A})$. In such a case denote the absolute cotangent complex of $F$ by $L_F$. This a $\mathcal{O}_F$-module. $L_F$ is an object in the stable $\infty$-category $\hbox{lim}_{\hbox{Spec}A \rightarrow F} \hbox{Stab}(\mathcal{C}^{op}_{/A})$.

Suppose there is a map of $\infty$-prestacks $F \rightarrow F'$. Since $A\oplus N \rightarrow A$ has a canonical section, given by the zero derivative, $\hbox{Def}_F(x,N)$ is a pointed set. Denote by $$\hbox{Def}_{F/F'}(x,-): \hbox{Mod}_A \rightarrow \hbox{SSets}$$ the homotopy fiber of the map

$$df: \hbox{Def}_F(x,-) \rightarrow \hbox{Def}_{F'}(x,-)$$

There is an alternate description of $\hbox{Def}_{F/F'}(x,-): \hbox{Mod}_A \rightarrow \hbox{SSets}$. Consider the functor $G: \mathcal{C}\hbox{-stacks}_{/F'} \rightarrow \hbox{SSets}$ which is the restriction of $F$ to along the natural map $\mathcal{C}\hbox{-stacks}_{/F'} \rightarrow \mathcal{C}\hbox{-stacks}$. Then for a point $x: \hbox{Spec}A \rightarrow F$, there is a point $x:\hbox{Spec}A \rightarrow G$ where $\hbox{Spec}A$ is considered an object in the over-category $\mathcal{C}\hbox{-stacks}$ via the map $\hbox{Spec}A \rightarrow F \rightarrow F'$. The relative deformation functor at $x$, $\hbox{Def}_{F/F'}(x,-)$ is then equivalent toi the absolute deformation functor $\hbox{Def}_G(x,-)$.

$F \rightarrow F'$ has a relative cotangent complex at $x$ if $\hbox{Def}_{F/F'}(x,-)$ is corepresentable by an $n$-connective $A$-module $L_{F/F',x}$ for some integer $n$. 

$F \rightarrow F'$ has a {\it relative cotangent complex} if $F \rightarrow F'$ has a relative cotangent complex at $x$ for all points $x$ and given a commutative diagram in $\infty\hbox{-stacks}$

$$\xymatrix{
\hbox{Spec}A \ar[dr]_x \ar[rr]^u &&\hbox{Spec}B \ar[dl]^{x'}\\
&F
}$$

the natural morphism $u^*L_{F/F',x'} \rightarrow L_{F/F',x}$ is an equivalence in $\hbox{Mod}_A$.

Suppose there is a sequence of maps of $\infty$-prestacks
$$F \rightarrow F' \rightarrow F"$$
and suppose the relative cotangent complex $F'/F"$ exists, then there is an exact triangle in the stable $\infty$-category of $F$-modules 
$$L_{F'/F"}|F \rightarrow L_{F/F"} \rightarrow L_{F/F'}$$ in the sense that if either of the second or the third term exist then so does the other and the triangle.

\section{Obstruction Theory}

In this section we extend the To\"en-Vessozi \cite{HAG} {\it obstruction theory} formalism for derived affine schemes to algebraic $\infty$-stacks.

Suppose $d_{\eta}:X[M] \rightarrow X$ is a derivation, induced by a map $\eta:L_A \rightarrow M$ in $\hbox{Stab}(\mathcal{C}^{op}_{/A})$. Define $X_{\eta}[\Omega M] := \hbox{Spec}(A \oplus_{\eta} \Omega M)$. Then the pullback square 

$$\xymatrix{
A\oplus_{\eta}\Omega M \ar[r] \ar[d] &A \ar[d]^{d_0} \\
A \ar[r]^{d_{\eta}} &A\oplus M
}$$

means $X_{\eta}[\Omega M]$ is the homotopy pushout $X\coprod^h_{X[\Omega M]} X$ in the $\infty$-category of affine $\mathcal{C}$-schemes.  

\begin{definition}(\cite{HAG})
An $\infty$-prestack $F$ {\it has an obstruction theory} if 
\begin{enumerate}[(i)]
\item $F$ is infinitesimally cohesive
\item $F$ has a cotangent complex
\end{enumerate}
\end{definition}

Geometric $\infty$-stacks always have an obstruction theory.

Suppose $F$ is has an obstruction theory then there exists a natural obstruction $\alpha(x) \in \hbox{Hom}_{\hbox{Mod}_A}(L_{F,x}, M)$ for a $A$-point $x:X\rightarrow F$ and $X_{\eta}[\Omega M]$ as defined above. This cohomological (Andre-Quillen) obstruction vanishes iff the dotted arrow exists in the diagram

$$\xymatrix{
X_{\eta}[\Omega M] \ar@{-->}[dr]^{x'}\\
X \ar[u] \ar[r]_x &F 
}$$

If $\alpha(x)=0$, the space of lifts of $x$, $\hbox{Hom}_{X/Aff_{\mathcal{C}}}(X_{\eta}[\Omega M], F)$, is isomorphic to $\hbox{Hom}_{\hbox{Mod}_A}(L_{F,x}, \Omega M) \simeq$ $$\Omega \hbox{Hom}_{\hbox{Mod}_A}(L_{F,x},M).$$

Is there a similar obstruction theory for lifting a family of object over an algebraic $\infty$-stack classified by a moduli stack $F$ which has an obstruction theory? Suppose $X = \hbox{colim}\hbox{Spec}A_{\bullet}$ (colimit in the $\infty$-category $\mathcal{C}$, i.e. the category of affine $\mathcal{C}$-schemes) where $A_{\bullet}$ is cosimplicial $\mathcal{C}^{op}$-object. Let $N \in \hbox{Stab}(\mathcal{C}^{op}_{/A_{\bullet}}$ be a $A_{\bullet}$-module and let $A^{\eta}_{\bullet}$ be the square-zero extension of $A_{\bullet}$ along a derivation $\eta:L_X \rightarrow N$. We want to find an obstruction for existence of the dotted arrow in 

$$\xymatrix{
\hbox{Spec}(A^{\eta}_{\bullet}) \ar@{-->}[dr]^{x'}\\
\hbox{Spec}(A_{\bullet}) \ar[u] \ar[r]_{x} &F
}$$

where $x$ is a $X$-point of $F$. It is clear from definitions that $\hbox{Spec}(A^{\eta}_{\bullet}) \simeq \hbox{Spec}A_{\bullet} \coprod^h _{\hbox{Spec}A_{\bullet}\oplus N} \hbox{Spec}A_{\bullet}$.

We need to verify that the following is an equivalence of simplicial sets when $F$ is infinitesimally cohesive
$$F(A^{\eta}_{\bullet}) \simeq F(A_{\bullet}) \times^h_{F(A_{\bullet}\oplus N)} F(A_{\bullet}).$$ 

Here for any cosimplicial $\mathcal{C}$-object $B_{\bullet}$, $F(B_{\bullet})$ is defined to be $F(\hbox{colim}_{\Delta^{op}}\hbox{Spec}B_{\bullet})$ using the Kan extension along the Yoneda map $\mathcal{C} \rightarrow P(\mathcal{C})$.

The following sequence of equivalences gives our desired equivalence.

$$F(A^{\eta}_{\bullet}) \simeq \hbox{Tot}_{[n]\in \Delta} F(A^{\eta}_{[n]})$$
$$\simeq \hbox{Tot}_{[n]\in \Delta}(F(A_{[n]}) \times ^h _{F(A_{[n]} \oplus N_{[n]})} F(A_{[n]}))$$
$$\simeq \hbox{Tot} F(A_{\bullet}) \times ^h _{\hbox{Tot}F(A_{\bullet}\oplus N)} \hbox{Tot}F(A_{\bullet})$$
$$\simeq F(\hbox{colim}\hbox{Spec}A_{\bullet}) \times^h_{F(\hbox{colim}\hbox{Spec}(A_{\bullet}\oplus N))} F(\hbox{colim}\hbox{Spec}A_{\bullet})$$

\section{Moduli of compact objects of $QC(X)$}

\begin{definition}(\cite{HTT})
A object $x$ in an $\infty$-category $\mathcal{D}$ is {\it compact} if the functor $\hbox{Hom}_{\mathcal{D}}(x,-):\mathcal{D} \rightarrow \infty\hbox{-groupoids}$ commutes with small colimits; $$\hbox{Hom}_{\mathcal{D}}(x, \hbox{colim}_{\alpha}y_{\alpha}) \simeq \hbox{colim}_{\alpha}\hbox{Hom}_{\mathcal{D}}(x,y_{\alpha}).$$
$\mathcal{D}$ is {\it compactly generated} if there exists a family of compact objects $\{x_{\alpha} \}_{\alpha}$ such that, any map $X \rightarrow Y$ in $\mathcal{D}$ is an equivalence if and only if $\hbox{Hom}_{\mathcal{D}}(x_{\alpha}, X) \rightarrow \hbox{Hom}_{\mathcal{D}}(x_{\alpha}, Y)$ is an weak equivalence of simplicial sets for all $\alpha$.
\end{definition}

A stable $\infty$-category $\mathcal{D}$ is compactly generated if there is a family of compact objects such that $y \in \mathcal{D}$ is the zero object iff $\hbox{Hom}_{\mathcal{D}}(x_{\alpha}, y)$ is a contractible simplicial set for all ${\alpha}$. In other words, for any arbitrary $y$ which is not the zero object, there is a non-zero map $c \rightarrow y$ from some compact object $c$.


The $\infty$-stack of perfect quasi-coherent modules $QC^{perf}$. Consider the $\infty$ functor considered as an object in $P(\mathcal{C})$ $$\hbox{Mod}: \mathcal{C}^{op} \rightarrow \widehat{Cat}_{\infty,st}$$ $$A\mapsto \hbox{Stab}(\mathcal{C}^{op}_{/A})^{\omega}.$$

For a map $A \rightarrow B$ in $\mathcal{C}$, there is a map of $\infty$ categories $\hbox{Stab}(\mathcal{C}^{op}_{/A})^{\omega} \rightarrow \hbox{Stab}(\mathcal{C}_{/B})^{\omega}$ since compact object objects map to compact objects. This extends to an $\infty$-functor. $QC^{perf}$ is the $\infty$-stack (fppf topology over connective $E_{\infty}$ rings)

$$QC^{perf}: \infty-\hbox{stacks} \rightarrow \hat{Cat}_{\infty,st}$$ obtained by Kan extension along the Yoneda embedding.

The objects of $\hbox{Stab}(\mathcal{C}^{op}_{/A})^{\omega}$ will be called {\it perfect complexes} of modules over $A$.

The stack $QC^{perf}$ is key to understanding the question of compact generation of the stable $\infty$-category $QC(X)$. We need that $QC^{perf}$ has an obstrction theory. In order for this we need to establish two things about $QC^{perf}$

\begin{itemize}
\item $QC^{perf}$ is infinitesimally cohesive
\item $QC^{perf}$ has a cotangent complex
\end{itemize}

It follows from a result of To\"en-Vessozi \cite{HAG} that it is enough to show that
\begin{itemize}
\item $QC^{perf}$ is infintesimally cohesive
\item The diagonal $\Delta:QC^{perf} \rightarrow QC^{perf}\times QC^{perf}$ is $n$-geometric for some $n$.
\end{itemize}

The first follows from the fact that $QC$ is infinitesimally cohesive. For the second part, let $A \in \mathcal{C}^{op}$ and let $x,y$ be objects in $QC^{perf}(\hbox{Spec}A)$. In other words $x$ and $y$ are perfect modules over $A$. Let $\Omega_{x,y}QC^{perf}$ be the pullback in the $\infty$ category of $\infty$-stacks.

$$\xymatrix{
\Omega_{x,y}QC^{\omega} \ar[r] \ar[d]  &QC^{\omega} \ar[d]^{\Delta}\\
\hbox{Spec}A \ar[r]_{x,y} &QC^{\omega} \times QC^{\omega} 
}$$ 

We'll show that $\Omega_{x,y}QC^{\omega}$ is an algebraic $n$-stack ($n$-truncated) for some $n$ depending on $A$, $x$ and $y$. The proof is based on the Artin-Lurie criterion. 

\begin{theorem}(Lurie)
A functor $F: conn E_{\infty}-rings \rightarrow SSets$ is a derived algebraic $n$-stack (in Lurie's sense, $n$-truncated) iff the following are satisfied
\begin{enumerate}[(i)]
\item $F$ is a sheaf in the etale topology
\item $F$ is $\omega$-accessible, it preserves $\omega$-filtered colimits
\item $F$ is nilcomplete, carries Postnikov towers to limits
\item $F$ is infinitesimally cohesive
\item $F$ has a cotangent complex
\item $F$ is formally effective
\item The restriction of $F$ to discrete commutative rings factors through $\hbox{SSets}^{\leq n}$.
\end{enumerate}
\end{theorem}

We'll show the existence of the cotangent complex for $\Omega_{x,y}QC^{\omega}$. Checking the other hypotheses in the Artin-Lurie criterion are easy.

Let $B$ be an object under $A$ in $\mathcal{C}^{op}$. Then the restriction of the functor $\Omega_{x,y}QC^{\omega}$ to $\mathcal{C}^{op}_{A/}$ can be described as 
$$\Omega_{x,y}QC^{\omega} = \hbox{Map}_{\hbox{Mod}_B}(x\otimes_AB, y\otimes_AB).$$

Use the notation $\mathcal{F} = \Omega_{x,y}QC^{\omega}_{/\hbox{Spec}A}: \mathcal{C}^{op}_{A/} \rightarrow \hbox{SSet}$ for the restriction of the functor $\Omega_{x,y}QC^{\omega}: \mathcal{C}^{op} \rightarrow \hbox{SSet}$ along the natural functor $\mathcal{C}^{op}_{A/} \rightarrow \mathcal{C}^{op}$. The structure morphism $\hbox{Spec}A \rightarrow \hbox{Spec}S$ has a cotangent complex. Therefore in order to show that $\Omega_{x,y}QC^{\omega}$ has an absolute cotangent complex it is sufficient to show that $\Omega_{x,y}QC^{\omega} \rightarrow \hbox{Spec}A$ has a relative cotangent complex, which is simply the cotangent complex of $F$. 

Let $B \in \mathcal{C}^{op}_{A/}$, an object in $\mathcal{C}\hbox{-stacks}_{/\hbox{Spec}A}$. Let $z: \hbox{Spec}B \rightarrow \mathcal{F}$ a map in $\mathcal{C}\hbox{-stacks}_{/\hbox{Spec}A}$. We want to show that the functor $\hbox{Def}_{\Omega_{x,y}QC^{\omega}/A}(x,-): \hbox{Mod}_B \rightarrow \hbox{SSet}$ is corepresentable. Recall this is equivalent to the functor $\hbox{Def}_{\mathcal{F}}(x,-)$. 
Let $B\oplus M$ be the trivial square-zero extension of $B$ along $M \in \hbox{Mod}_B$. We have 

$$\mathcal{F}(\hbox{Spec}B) = \hbox{Map}_{\hbox{Mod}_B}(x \otimes_AB, y \otimes_AB) \simeq \hbox{Hom}_{\hbox{Mod}_A}(x,y\otimes_AB)$$

$$\mathcal{F}(\hbox{Spec}(B\oplus M)) = \hbox{Map}_{\hbox{Mod}_{(B\oplus M)}}(x\otimes_A(B\oplus M), y \otimes_A (B\oplus M)) \simeq \hbox{Hom}_{\hbox{Mod}_A}(x,y\otimes_A(B\oplus M))$$

$$\simeq \hbox{Hom}_{\hbox{Mod}_A}(x,y\otimes_AB) \times \hbox{Hom}_{\hbox{Mod}_A}(x,y\otimes_AM)$$

All these equivalences commute with the natural map $$\mathcal{F}(\hbox{Spec}(B\oplus M)) \rightarrow \mathcal{F}(\hbox{Spec}B).$$

Therefore the deformation space $\hbox{Def}_{\mathcal{F}}(x,M)$ which is the homotopy fiber of this map at $x$ is equivalent to $\hbox{Hom}_{\hbox{Map}_A}(x,y\otimes_A M)$. There is a chain of equivalences
$$\hbox{Def}_{\mathcal{F}}(x,M) \simeq \hbox{Hom}_{\hbox{Mod}_A}(x,y \otimes_AM)$$
$$\simeq \Omega^{\infty}(\hbox{Mor}_A(x,y)\otimes_A M)$$
$$\simeq \Omega^{\infty}((\hbox{Mor}_A(x,y)\otimes_AB)\otimes_BM))$$
$$\simeq \Omega^{\infty}((\hbox{Mor}_B((\hbox{Mor}_A(x,y)\otimes_AB)^v,M))$$
$$\simeq \hbox{Hom}_{\hbox{Mod}_B}(\hbox{Mor}_A(x,y)\otimes_AB)^v,M)$$

The notation $\hbox{Mor}_A(x,y)$ is used for $\hbox{Hom}_{\hbox{Mod}_A}(x,y)$ when considered as an object of the stable $\infty$-category $\hbox{Mod}_A$. 

The equivalences follow from the facts that $\hbox{Mor}_A(x,y)$ is a compact object when $x$ and $y$ are compact, $\hbox{Mod}_A$ is compactly generated under filtered colimits by $A$ and compact objects are dualizable in $\hbox{Mod}_B$.

Therefore $\hbox{Def}_{\mathcal{F}}(x,-)$ is corepresentable by the $B$-module $L_{\mathcal{F},x}:=(\hbox{Mor}_A(x,y)\otimes_AB)^v$.

Suppose given a commutative diagram in $\mathcal{C}\hbox{-stacks}_{/\hbox{Spec}A}$
$$\xymatrix{
\hbox{Spec}C \ar[rr]^u \ar[dr]_w &&\hbox{Spec}B \ar[dl]^z \\
&\mathcal{F}
}$$

we have the equivalences

$$L_{\mathcal{F},w} \simeq (\hbox{Mor}_A(x,y)\otimes_AC)^v \simeq \hbox{Mor}_A(\hbox{Mor}_A(x,y),C)$$
$$u^*L_{\mathcal{F},z} \simeq (\hbox{Mor}_A(x,y)\otimes_AB)^v\otimes C \simeq \hbox{Mor}_A(\hbox{Mor}_A(x,y),B)\otimes_CB$$

The equivalences follow simply from adjunction are compatible with the natural map $u^*L_{\mathcal{F},z} \rightarrow L_{\mathcal{F},w}$ making it an equivalence in $\hbox{Mod}_C$.

This completes the proof that $\Omega_{x,y}QC^{\omega}$ has a cotangent complex. We need to verify the rest of the Artin-Lurie conditions to show that it is an algebraic stack. Then applying the proposition of \cite{HAG} it follows that $QC^{\omega}$ has a cotangent complex.

\section{Proof of the Main Theorem}

\begin{definition}(\cite{Drinfeld})
A derived $\infty$-$\mathcal{C}$-stack $X$ is {\it perfect} if 
\begin{enumerate}[(i)]
\item $X$ has affine diagonal,
\item $QC(X)$ is a presentable stable $\infty$-category, or equivalenty the triangulated category $\hbox{ho}(QC(X))$ is compactly generated.
\end{enumerate}
\end{definition}

Suppose $A_{\bullet}$ is a cosimplicial object in $\mathcal{C}^{op}$ which is level-wise truncated as objects in the $\infty$-category $\mathcal{C}^{op}$. Then the derived algebraic $\mathcal{C}$-stack $X = \hbox{colim}_{\Delta^{op}} \hbox{Spec}(A_{\bullet})$ can be obtained as finitely many square-zero extensions of the (non-derived) {\it classical} algebraic $\infty$-$\mathcal{C}$-stack $$X^{cl} = \hbox{colim}_{\Delta^{op}}(\hbox{Spec}(\pi_0A_{\bullet}))$$

There is a natural map $i:X^{cl} \rightarrow X$. Suppose we know that $X^{cl}$ is perfect, what can be said about the perfectness of derived counterpart $X$? Since $X^{cl} \rightarrow X$ is an infintesimal extension of stacks, we shall consider the following question: suppose $i: X \rightarrow \widetilde{X}$ is a square-zero extension of an $\infty$-algebraic stack $X$ and suppose $QC(X)$ is compactly generated. What can be said about the presentability of the stable category $QC(\widetilde{X})$?

\begin{enumerate}[(I)]

\item We've seen in the previous section that $QC^{\omega}$ has an obstruction theory. Therefore we can use $L_{QC^{\omega}}$ to lift the compact objects in $QC(X)$ to compact objects in $QC(\widetilde{X})$. The space of all such lifts is a deformation space 

$$\xymatrix{
&\widetilde{X} \ar@{-->}[d]^{\tilde{u}}\\
X \ar[ur]^i \ar[r]_u &QC^{\omega}
}$$

There is an obstruction in the Andre-Quillen cohomology group $$\alpha(u) \in \hbox{Hom}_{\hbox{Mod}_{\mathcal{O}_X}}(u^*L_{QC^{\omega}}, N)$$ (where $N \in \hbox{Stab}(\mathcal{C}^{op}_{/A_{\bullet}})$ is a $\mathcal{O}_X$-module, so that $\widetilde{X} = \hbox{colim}\hbox{Spec}(A^{\eta}_{\bullet})$ for some derivation $\eta:L_X \rightarrow N$). If $\alpha(u)=0$ let $\tilde{u}$ be a deformation of $u$.

\item Given  $x \in QC(\widetilde{X})$. Then $i^*(x) \in QC(X)$. Since $QC(X)$ is compactly generated, there exists $u \in QC(X)^{\omega}$ and a non-zero map $u \rightarrow i^*(x)$ in $QC(X)$.  We want to know if there is a lift of the map $f: u \rightarrow i*(x)$ in $QC(X)$ to an map $\tilde{u} \rightarrow x$ in $QC(\widetilde{X})$ under the map of stable $\infty$-categories $$QC(\widetilde{X}) \rightarrow QC(X)$$ induced by the natural map $i:X \rightarrow \widetilde{X}$.

The space of all possible lifts is the space of deformations of the map $u \rightarrow i^*(x)$ and is controlled by the cotangent complex of the $\infty$-stack $\Omega_{u,i^*x}QC$.

We'll give an description of the space of lifts of the map $f:u \rightarrow i^*(x)$ to $\tilde{f}: \tilde{u} \rightarrow x$ in $QC(\widetilde{X})$.

That $\Omega_{u,i^*(x)}QC \simeq X \times_{QC} X$ in the category of $\infty$-stacks means that for any affine $\mathcal{C}$-scheme $\Omega_{u,i^*(x)}QC(\hbox{Spec}A)$ is the $\infty$ category $\hbox{Hom}_{\infty\hbox{-stacks}}(\hbox{Spec}A, \Omega_{u,i^*x}QC)$ in which the $0$-simplices are triplets $(f,g,\phi)$ where $$f,g:\hbox{Spec}A \rightarrow X$$ and $$\phi:f^*u \rightarrow g^*i^*(x)$$ is a map in $\hbox{Mod}_A$. The $1$-cells are morphisms between such triplets defined in the natural way.

In particular, the if we take the test space to be $X$ itself and a square zero-extension $\widetilde{X}$ of $X$, then the mapping spaces are 
$$ \Omega_{u,i^*x}QC(X) = \hbox{Hom}_{\infty\hbox{-St}}(X, X \times_{QC} X)$$
$$ \Omega_{u,i^*x}QC(\widetilde{X}) = \hbox{Hom}_{\infty\hbox{-St}}(\widetilde{X},X \times_{QC} X)$$

The first space is the $\infty$-category whose objects are triplets $(f,g:X\rightarrow X, \phi: f^*x \rightarrow g^*y \in \hbox{Mod}_{\mathcal{O}_X})$. The second space is the $\infty$-category whose objects are triplets $(f',g':\widetilde{X}\rightarrow X, \phi ': f'^*x \rightarrow g'^*y \in \hbox{Mod}_{\mathcal{O}_{\widetilde{X}}})$. Here $f'^*x$, $g'^*y$ and $\phi '$ are {\it not} deformations of $f^*x$, $g^*y$ and $\phi$ respectively. However if we consider the point in $\Omega_{u, i^*x}QC(X)$ represented by the object $(1,1,f)$ corresponds to the triplet $(x,y,f:u \rightarrow i^*x)$ in $X \times_{QC}X$, then the fiber of $\Omega_{u,i^*x}QC(\widetilde{X}) \rightarrow \Omega_{u,i^*x}QC(X)$ over this point

$$\xymatrix{
\mathfrak{D} \ar[r] \ar[d] &\Omega_{u,i^*x}QC(\widetilde{X}) \ar[d] \\
\star \ar[r]_{(1,1,f)} &\Omega_{u,i^*x}QC(X)
}$$

is the $\infty$-category of objects $(u',x',\tilde{f})$ which are respectively deformations of $x$, $u$ and $f:u\rightarrow i^*(x)$ to $QC(\widetilde{X})$. This deformation space is larger than the one we need. We want the space of deformations of the map $f$ that keeps a fixed choice of deformations of the source $u$ and target $i^*x$.

Consider the moduli functor $\mathcal{F}:\infty\hbox{-stacks}_{/X \times X} \rightarrow {Cat}_{\infty}$ obtained by restricting $\Omega_{u,i^*x}QC$ along the natural functor $\infty\hbox{-stacks}_{/X\times X} \rightarrow \infty\hbox{-stacks}$.

Let $z: \hbox{Spec}A \rightarrow \mathcal{F}$ be a map in $\infty\hbox{-St}_{/X \times X}$. Then the mapping space
$$\mathcal{F}(\hbox{Spec}A)$$ is the $\infty$-category whose objects are maps $\phi:f^*x \rightarrow g^*y$ in $\hbox{Mod}_{\mathcal{O}_X}$. Here $f,g:\hbox{Spec}A \rightarrow X \times X$ is the test space in $\infty\hbox{-St}_{X \times X}$. Denote this test space by $\hbox{Spec}A_{f,g}$. 

$X$ is naturally an object in $\infty\hbox{-St}_{/X\times X}$ via the identity maps. We will denote this version of $X \in \infty\hbox{-St}_{/X \times X}$ by $X_{1,1}$. 

Let $\widetilde{X}$ be considered an object in $\infty\hbox{-St}_{/X\times X}$ via the derivations $d_1,d_2:\widetilde{X} \rightarrow X \times X$ so that $d_1^*u = \tilde{u}$ and $d_2^*(i^*x) = x$. Denote this object of the over category by $\widetilde{X}_{d_1,d_2}$.

Now consider the point in $\mathcal{F}(X_{1,1})$ corresponding to the map $f:u \rightarrow i^*x$. The fiber of the natural map $\mathcal{F}(\widetilde{X}_{d_1,d_2}) \rightarrow \mathcal{F}(X_{1,1})$ over this point

$$\xymatrix{
\mathfrak{D}_0 \ar[r] \ar[d] &\mathcal{F}(\widetilde{X}_{d_1,d_2}) \ar[d]\\
\star \ar[r]_f  &\mathcal{F}(X_{1,1})
}$$

is the $\infty$-category whose objects are exactly the deformations of the map $f:u \rightarrow i^*x$ in $QC(X)$ to $\tilde{f}: \tilde{u} \rightarrow x$ in $QC(\widetilde{X})$.

Therefore if $\mathcal{F}$ has an obstruction theory, this deformation problem of lifting the map $f$

$$\xymatrix{
&\widetilde{X}_{d_1,d_2} \ar@{-->}[d]^{\tilde{f}} \\
X_{1,1} \ar[ur]^i \ar[r]_f &\mathcal{F}
}$$

is controlled by the cotangent complex $L_{\mathcal{F}}$. Recall that this equivalent to the relative cotangent complex $L_{\Omega_{u,i^*x}QC/X\times X}$ with respect to the natural map $\Omega_{u,i^*x}QC = X \times_{QC} X \rightarrow X \times X$ of $\infty$-stacks.

More precisely, there is a cohomological obstruction $$\beta(f) \in \hbox{Hom}_{\mathcal{O}_{X_{1,1}}}(f^*L_{\mathcal{F}}, N)$$ Alternately this obstruction lives in $$\hbox{Hom}_{\mathcal{O}_X}((1,1,f)^*L_{\Omega_{u,i^*x}QC},N).$$

If $\beta(f)=0$ there exists deformations of $f$. The space of all possible deformations $\tilde{f}:\tilde{u} \rightarrow x$ is $$\Omega\hbox{Hom}_{\mathcal{O}_X}((1,1,f)^*L_{\Omega_{u,i^*x}QC},N).$$

\end{enumerate}

For these two steps to work we need the two moduli stacks $QC^{\omega}$ and $\Omega_{u,i^*x}QC$ {\it have deformation theory}. In other words that they are infinitesimally cohesive and have cotangent complexes. This has already established for $QC^{\omega}$. Checking that the second space is infintesimally cohesive is formal. Now we come to the existence of the cotangent complex for $\Omega_{u,i^*x}QC$. 

Since $i^*x$ need not be {\it compact}, $\Omega_{u,i^*x}QC$ {\it does not} have a cotangent complex in general. 

However $i^*x \in \hbox{Mod}_{\mathcal{O}_X}$ and $X$ is perfect. Therefore $i^*x$ is a filtered colimit of perfect modules over $\mathcal{O}_X$. Let us suppose that $i^*x = \hbox{colim} y_{\beta}$, for $\beta: X \rightarrow QC^{\omega}$. Then the natural map

$$\hbox{Hom}_{\mathcal{O}_X}(u, \hbox{colim} y_{\alpha}) \rightarrow \hbox{colim}\hbox{Hom}_{\mathcal{O}_X}(u,y_{\alpha})$$ is an equivalence since $u$ is compact. Therefore any map $f:u \rightarrow \hbox{colim}y_{\alpha}= i^*x$ factors through $\beta: u \rightarrow y_{\beta}$ for some $\beta$.

Since $d_2^*$ is a left adjoint, it preserves colimits, $$x = d_2^*(\hbox{colim}y_{\alpha}) \simeq \hbox{colim}(d_2^*y_{\alpha})$$
It is clear that $d_2^*y_{\alpha}$ need not be compact. 


Replace the moduli stacks $\Omega_{u, i^*x}QC$ in the second step with $\Omega_{u,y_{\beta}}QC$. This one does indeed have an obstruction theory. This means that the functor 

$$\mathcal{G}: \mathcal{C}_{/X \times X} \rightarrow \widehat{Cat}_{\infty}$$ is infinitesimally cohesive and has a cotangent complex.

There exists a natural obstruction in the Andre-Quillnen cohomology $$\alpha(u,y_{\beta}) \in \hbox{Hom}_{\mathcal{O}_X}(\beta^*L_{\mathcal{G}},N)$$ for lifting the map $\beta:u \rightarrow y_{\beta}$ to $\widetilde{\beta}:\widetilde{u} \rightarrow d_2^*(y_{\beta})$. The space of all such deformations is equivalent to the space 

$$\Omega\hbox{Hom}_{\mathcal{O}_X}(\beta^*L_{\mathcal{G}},N)$$ with loops based at the trivial derivation.

$i = \hbox{colim}d_2^*(y_{\beta})$ implies there is a unique map $d_2^*(y_{\beta}) \rightarrow x$. Compose this with $\widetilde{\beta}$ to obtain the desired lift of $u \rightarrow i^*x$ to $\widetilde{X}$.


\addcontentsline{toc}{section}{References}

\begin{thebibliography}{99}

\bibitem{Drinfeld} {\bf David Ben-Zvi, John Francis, David Nadler},
{\em Integral Transforms and Drinfeld Centers in Derived Algebraic Geometry}
(arXiv:0805.0157)

\bibitem{DAGIV} {\bf Jacob Lurie},
{\em Derived Algebraic Geometry IV}

\bibitem{HTT} {\bf Jacob Lurie},
{\em Higher Topos Theory}, Annals of Mathematics Studies.
(arXiv:math/0608040)

\bibitem{Neeman} {\bf A. Neeman},
{\em The Grothendieck duality via Bousefield's technique's and Brown representability}
(arXiv:alg-geom/9412022)

\bibitem{TT} {\bf A. Neeman},
{\em The connection between the K-theory localization theorem of Thomason,
Trobaugh and Yao and the smashing subcategories of Bousfield and Ravenel}.
Annales scientifiques de l’É.N.S. 4e série, tome 25, no 5 (1992)


\bibitem{HAG} {\bf B. To\"en, G. Vezzosi},
{\em Homotopical Algebraic Geometry II}.
(arXiv:math/0404373)


\end{thebibliography}
\end{document}